\documentclass[11pt]{amsart}
\usepackage{times,epsf,amssymb,amsmath}

\begin{document}

\newtheorem{thm}{Theorem}[section]
\newtheorem{lem}[thm]{Lemma}
\newtheorem{cor}[thm]{Corollary}

\theoremstyle{definition}
\newtheorem{defn}{Definition}[section]

\theoremstyle{remark}
\newtheorem{rmk}{Remark}[section]

\def\square{\hfill${\vcenter{\vbox{\hrule height.4pt \hbox{\vrule
width.4pt height7pt \kern7pt \vrule width.4pt} \hrule height.4pt}}}$}

\def\T{\mathcal T}

\newenvironment{pf}{{\it Proof:}\quad}{\square \vskip 12pt}

\title{Least Area Planes in Hyperbolic $3$-Space are Properly Embedded}
\author{Baris Coskunuzer}
\address{Koc University \\ Department of Mathematics \\ Sariyer, Istanbul 34450 Turkey}
\email{bcoskunuzer@ku.edu.tr}
\thanks{The author is partially supported by NSF Grant DMS-0603532}

\maketitle


\newcommand{\SH}{S^2_{\infty}(\mathbf{H}^3)}
\newcommand{\PI}{\partial_{\infty}}
\newcommand{\SI}{S^2_{\infty}}
\newcommand{\BHH}{\mathbf{H}^3}
\newcommand{\CH}{\mathcal{C}(\Gamma)}
\newcommand{\BH}{\mathbf{H}}
\newcommand{\BR}{\mathbf{R}}
\newcommand{\BC}{\mathbf{C}}
\newcommand{\BZ}{\mathbf{Z}}

\begin{abstract}

We show that if $\Sigma$ is an embedded least area (area minimizing) plane in $\BH^3$ whose
asymptotic boundary is a simple closed curve with at least one smooth point, then $\Sigma$ is
properly embedded in $\BH^3$.

\end{abstract}

\section{Introduction}

The asymptotic Plateau problem in hyperbolic space asks the existence of a least area (area minimizing) plane $\Sigma \subset \BHH$
asymptotic to given simple closed curve $\Gamma\subset \SH$. This problem is solved by Michael Anderson in his seminal papers \cite{A1} ,
\cite{A2}. He proved the existence of a solution for any given simple closed  curve in the sphere at infinity. Later, by using topological
techniques, Gabai proved a similar result for $\BHH$ with any cocompact metric in \cite{Ga}. Then, the author generalized these results to
Gromov hyperbolic $3$-spaces with cocompact metric \cite{Co1}.

Properly embeddedness of the solution has been questioned by both Anderson and Gabai. Nevertheless, only known results about the properly
embeddedness of least area planes in $\BHH$ is the existence of some properly embedded least area plane for a given simple closed curve in
$\SH$ by \cite{So1}, \cite{So2}, \cite{Co3}. It is still not known if there exists a nonproperly embedded least area plane in $\BHH$ whose
asymptotic boundary is a simple closed curve in $\SH$.

On the other hand, recently Colding and Minicozzi proved a very powerful result about properly
embeddedness of complete embedded minimal disks in $\BR^3$ in \cite{CM}. They proved Calabi-Yau
Conjectures for embedded surfaces by relating intrinsic distances and extrinsic distances of the
minimal disk. As a corollary, they proved that any complete embedded minimal plane in $\BR^3$ must
be proper.

In this paper, we prove an analogous result in $\BH^3$. We show that if $\Sigma$ is an embedded
least area plane in $\BHH$ whose asymptotic boundary is a simple closed curve with at least one
smooth point, then $\Sigma$ is properly embedded in $\BH^3$. Instead of relating the intrinsic and
extrinsic distances as Colding and Minicozzi did, we use powerful topological arguments. The main
result of the paper is as follows:

\vspace{0.3cm}

\noindent \textbf{Theorem 5.2.} Let $\Sigma$ be a complete embedded least area plane in $\BHH$ with
$\PI\Sigma=\Gamma$ where $\Gamma$ is a simple closed curve in $\SH$ with at least one smooth
($C^1$) point. Then, $\Sigma$ must be proper.

\vspace{0.3cm}

The organization of the paper is as follows: In the next section we will cover some basic results which will be used in the following
sections. In section 3, we will analyze the intersection of a least area plane $\Sigma$ with balls exhausting $\BHH$. Then in section 4, we
will prove the key lemma which is the most important step for the main result. In section 5, we will prove the main result. Finally in
section 6, we will have some concluding remarks.

\subsection{Acknowledgements:}

I would like to thank David Gabai and Yair Minsky for very useful conversations.

\section{Preliminaries}

In this section, we will overview the basic results which we will use in the following sections. First, we will give the definitions of
least area (area minimizing) planes.

\begin{defn}A {\em least area disk} is a disk which has the smallest area among the disks with the same
boundary. A {\em least area plane} is a plane such that any compact subdisk in the plane is a least
area disk.
\end{defn}

\begin{defn} An immersed surface $S$ in $\BHH$ is {\em proper} if the preimage of any compact subset of
$\BHH$ is compact in the surface $S$. If an embedded surface $S$ in $\BHH$ is proper, we will call
$S$ as {\em properly embedded}.
\end{defn}

\begin{defn} Let $A$ be a subset of $\SH$. Then the \textit{convex hull} of $A$, $CH(A)$, is the
smallest closed convex subset of $\BHH$ which is asymptotic to $A$. Equivalently, $CH(A)$ can be defined as the intersection of all
supporting closed half-spaces of $\BHH$ \cite{EM}.
\end{defn}

It is a well-known fact in minimal surface theory that if $M$ is a minimal surface in $\BHH$ with
$\PI M= \Gamma$, then $M\subset CH(\Gamma)$.

Now, we will quote the basic results on asymptotic Plateau problem.

\begin{thm} \cite{A2} Let $\Gamma$ be a simple closed curve in $\SH$. Then there exist
a complete least area plane $\Sigma$ in $\BHH$ asymptotic to $\Gamma$ at infinity.
\end{thm}

Later, Hardt and Lin showed regularity at infinity for these solutions in \cite{HL}. Then, in
\cite{To}, Tonegawa generalized this result to any complete constant mean curvature hypersurfaces
in $\BH^n$. The following theorem is indeed true for any dimension. For simplicity, we only mention
the result which we are interested in.

\begin{thm} \cite{HL}, \cite{To}
Let $\Gamma$ be a $C^1$ regular simple closed curve in $\SH$, and $\Sigma$ be a least area plane in
$\BHH$ asymptotic to $\Gamma$. Let $\overline{\BHH}$ be the compactification of $\BHH$ with $\SH$.
Then, there is a neighborhood $X$ of $\SH$ in $\overline{\BHH}$ such that $(\Sigma\cup\Gamma)\cap
X$ is a finite union of $C^1$ submanifolds of $\overline{\BHH}$.
\end{thm}

The following is a simple lemma which states that the intersection of a least area plane with a
ball is generically a disjoint union of disks.

\begin{lem} \cite{A2}
Let $\Sigma$ be a complete minimal plane with $\PI \Sigma = \Gamma$ where $\Gamma$ is a simple
closed curve in $\SH$. Then for almost all $r>0$, $\Sigma \cap B_r(0)$ is a disjoint union of
disks.
\end{lem}

The following lemma will be used later. The proof basically uses Meeks-Yau exchange roundoff trick
for least area disks \cite{MY2}

\begin{lem} \cite{Co2}
Let $\Gamma_1$ and $\Gamma_2$ be two disjoint simple closed curves in $\SH$. If $\Sigma_1$ and
$\Sigma_2$ are least area planes in $\BH^3$ with $\PI \Sigma_i = \Gamma_i$, then $\Sigma_1$ and
$\Sigma_2$ are disjoint, too.
\end{lem}

\section{Intersection of Least Area Planes with Balls Exhausting $\BHH$}

In this section, we will analyze the intersection of a given least area plane $\Sigma$ with balls
with fixed center and increasing radius in $\BHH$. By Lemma 2.3, we know that for a generic radius
$r>0$, the intersection $B_r(0)\cap\Sigma$ is a collection of disjoint disks. If $\Sigma$ is not
proper, we will show that there is a $r_0$ such that for a generic $r>r_0$, $B_r(0)\cap\Sigma$
contains infinitely many disjoint disks. In this section, we will analyze these disks, and classify
them accordingly.

Let $\Sigma$ be a least area plane with $\PI \Sigma = \Gamma$ where $\Gamma$ is a simple closed
curve in $\SH$. Fix a point $0$ on $\Sigma \subset CH(\Gamma)$. Let $B_r(0)$ be a closed ball with
radius $r$ (extrinsic) and center $0$ in $\BHH$. By Lemma 2.3, we know that for a generic radius
$r>0$, the intersection $B_r(0)\cap\Sigma$ is a collection of disjoint disks.

\begin{lem}
Let $\Sigma$ be an embedded least area plane in $\BHH$ with $\PI \Sigma = \Gamma$ where $\Gamma$ is
a simple closed curve in $\SH$. If $\Sigma$ is not proper, then there exist $r_0>0$ such that for a
generic $r>r_0$, $B_r(0)\cap\Sigma$ contains infinitely many disjoint disks.
\end{lem}

\begin{pf}
Let $\varphi: D^2\rightarrow \BHH$ be the smooth embedding with $\varphi(D^2)= \Sigma$. If $\Sigma$
is not proper, then there exist a compact subset $K$ of $\BHH$ such that $\varphi^{-1}(K)=E$ is not
compact in $D^2$. This implies $E$ is not bounded in $D^2$.

Let $r_0>0$ be a generic radius with $K\subset B_{r_0}(0)$. If $\varphi^{-1}(B_{r_0}(0))=E_{r_0}$, then clearly $E\subset E_{r_0}$ in
$D^2$. By genericity, the intersection of $B_{r_0}(0)\cap \Sigma$ is a collection of disjoint disks. Assume that there are finitely many
disks in the intersection, i.e $B_{r_0}(0)\cap \Sigma = \bigcup_{i=1}^N D_i$ where $D_i$ is a closed disk in $\Sigma$. Let $C_i$ be the
boundary of the disk $D_i$, i.e. $C_i=\partial D_i$. Consider $\gamma_i=\varphi^{-1}(C_i)$ in $D^2$. Since $\gamma_i$ is a simple closed
curve in $D^2$, it will enclose a bounded disk in $D^2$ which is the preimage of $D_i$. Hence, if the disks in the intersection
$B_{r_0}(0)\cap \Sigma$ are finitely many, then the preimage $\varphi^{-1}(B_{r_0}(0))=E_{r_0}$ must be bounded in $D^2$. This is a
contradiction as $E$ is not bounded in $D^2$ and $E\subset E_{r_0}$.

Since $r_0$ is any generic radius with $K\subset B_{r_0}(0)$, for any generic $r>r_0$,
$B_r(0)\cap\Sigma$ contains infinitely many disjoint disks.
\end{pf}

Now, we will categorize these infinitely many disks in the intersection of nonproper least area plane $\Sigma$ and sufficiently large ball
$B_r(0)$ in $\BHH$. Consider the intersection $CH(\Gamma)\cap B_r(0)$ for sufficiently large $r$, which is a convex body in $\BHH$. Let
$\partial^+ CH(\Gamma)$ and $\partial^- CH(\Gamma)$ be the two connected components of the boundary of $CH(\Gamma)$. Let $A_r$ be the
annulus in $\partial B_r(0)$ such that $\partial A_r = \delta^+_r\cup \delta^-_r$ where $\delta^\pm_r \subset \partial^\pm CH(\Gamma)\cap
\partial B_r(0)$, and $int(CH(\Gamma)) \cap \partial B_r(0) \subset A_r$.

Since $\Sigma \subset CH(\Gamma)$, for any disk $D$ in the intersection $B_r(0)\cap \Sigma$,
$\partial D$ must belong to $A_r$. We call a disk $D$ in the intersection $B_r(0)\cap \Sigma$ {\em
separating} if $\partial D$ is essential in $A_r$, and {\em nonseparating} otherwise.

Next, we will analyze the separating and nonseparating disks in the intersection of a ball with the
nonproper least area plane. This analysis will play essential role in the main result.

\section{Key Lemma}

In this section, we will prove the key lemma which is the most important step for the main result.
Roughly, we will show that if a disk $D$ is nonseparating with $\partial D \subset A_r$, then $D$
stays close to $A_r$, and  it does not come near $0$.

\vspace{0.3cm}

\noindent \textbf{Key Lemma:} Let $\Sigma$ be an embedded least area plane in $\BHH$ with $\PI
\Sigma = \Gamma$ where $\Gamma$ is a simple closed curve with at least one smooth ($C^1$) point in
$\SH$. Let $D_r$ be a nonseparating disk in $B_r(0)\cap\Sigma$. Then there is a function $F$  which
is a monotone increasing function with $F(r)\rightarrow \infty$ as $r \rightarrow \infty$, such
that $d(0,D_r)>F(r)$ where $d$ is the distance.

\vspace{0.3cm}

\begin{pf} First, we will give an outline of the proof. Then, we will prove the lemma in $2$ steps.\\

\noindent \textbf{Outline:} The idea to show that the nonseparating disks $D_r$ cannot be too close to the center is to construct a
barrier. By using the smooth point assumption, we will show the existence of a complete least area annuli $\mathcal{A}_s$ "linking"
$\Gamma$ such that $\PI \mathcal{A}_s = \Gamma^+_s\cup \Gamma^-_s$ in $\SH$ with $\Gamma^\pm_s\rightarrow \Gamma$ as $s\rightarrow\infty$.
Then, such an $\mathcal{A}_s$ will separate $\BHH$ into two parts, where one part contains $A_r$ which contains $\partial D_r$, and other
part contains $0$. Then, since $\mathcal{A}_s$ is least area annulus, and $D_r$ is least area disk, by exchange roundoff trick \cite{MY2},
$D_r$ cannot intersect $\mathcal{A}_s$. So, $\mathcal{A}_s$ becomes a barrier between $D_r$ and $0$. Since $\mathcal{A}_s \rightarrow
\infty$ as $\PI\mathcal{A}_s=\Gamma\pm_s\rightarrow\Gamma$, this defines a monotone increasing function
$F(r)$ with $F(r)\rightarrow \infty$ as $r \rightarrow \infty$, such that $d(0,D_r)>F(r)$.\\

\noindent \textbf{Step 1:} ({\em Existence of Least Area Annuli}) Let $\Gamma^+$ and $\Gamma^-$ be
two simple closed curves in opposite sides of $\Gamma$ in $\SH$ and sufficiently close to $\Gamma$.
Then there exists a complete least area annulus $\mathcal{A}$ in $\BHH$ with $\PI \mathcal{A} =
\Gamma^+ \cup
\Gamma^-$.\\

\begin{pf}
Let $x$ be a $C^1$ smooth point in $\Gamma$. Then, by the discussion in Section $1$ in \cite{HL},
we can find round circles $\gamma^+$ and $\gamma^-$ in the opposite sides of $\Gamma$ so that
$\gamma^+$ and $\gamma^-$ are as close as we want to $x\in\Gamma$. Let $P^+$ and $P^-$ be totally
geodesics planes asymptotic to round circles $\gamma^+$ and $\gamma^-$, respectively. Since we can
make $\gamma^+$ and $\gamma^-$  as close as we want by using the construction of Hass in \cite{Ha},
we can find curves $\alpha^+$ and $\alpha^-$ on $P^+$ and $P^-$ so that they cobound an embedded
least area annulus $\Delta$ whose area is strictly less than the area of the two totally geodesic
disks bounded by $\alpha^\pm$ on $P^\pm$ \cite{MY1}.

Now, let $N(\Gamma)$ be a neighborhood of $\Gamma$ which is an annulus in $\SH$ so that $\gamma^+$ and $\gamma^-$ are disjoint from
$N(\Gamma)$. Let $N^+(\Gamma)$ and $N^-(\Gamma)$ be the components of $N(\Gamma) - \Gamma$ in $\SH$. Foliate $N^\pm(\Gamma)$ with $C^1$
pairwise disjoint simple closed curves $\{\Gamma^\pm_s\}$ where $s\in(C,\infty)$, and $\Gamma^\pm_s\rightarrow\Gamma$ as
$s\rightarrow\infty$. We claim that for any pair $\Gamma^+_s$ and $\Gamma^-_s$, there is a complete least area annulus $\mathcal{A}_s$ in
$\BHH$ with $\PI \mathcal{A}_s = \Gamma^+_s\cup\Gamma^-_s$.

First, we fix a $s\in(C,\infty)$. By Theorem 2.1, for any simple closed curve in $\SH$, there exist a complete least area plane in $\BHH$.
Let $\Sigma^+_s$ and $\Sigma^-_s$ be the least area planes in $\BHH$ with asymptotic boundary $\Gamma^+_s$ and $\Gamma^-_s$, respectively.
Since $\Gamma^\pm_s$ and $\gamma^\pm$ are disjoint in $\SH$, then by Lemma 2.4, the least area planes $\Sigma^\pm_s$ and the geodesic
planes $P^\pm$ are pairwise disjoint. Hence, the annulus $\Delta$ intersects $\Sigma^\pm_s$ transversely in simple closed curves. Let
$\Omega^+_{s0}$ be a sufficiently large disk in $\Sigma^+_s$ with $\Delta \cap \Sigma^+_s \subset \Omega^+_{s0}$. Let $\beta^+_{s0}=
\partial \Omega^+_{s0}$ be simple closed curve in $\Sigma^+_s$. Since $\Gamma^+_s$ is $C^1$ regular, by Theorem 2.2, $\Sigma^+_s$ behave
nicely near asymptotic boundary. We foliate $\Sigma^+_s - \Omega^+_{s0}$ with pairwise disjoint simple closed curves $\beta^+_{st}$ where
$t\in [0,\infty)$. Similarly, we define $\Omega^-_{s0}$ and $\beta^-_{s0}$ in $\Sigma^-_s$, and foliate $\Sigma^-_s - \Omega^-_{s0}$ with
pairwise disjoint simple closed curves $\beta^-_{st}$ where $t\in [0,\infty)$.

Now, we claim that each pair $\beta^+_{st}$ and $\beta^-_{st}$ cobounds a least area annulus in $\BHH$ for any $t$. To prove that, we need
to show that there is an annulus with boundary $\beta^+_{st}\cup\beta^-_{st}$ whose area is less than the sum of the areas of least area
disks bounded by $\beta^+_{st}$ and $\beta^-_{st}$, say $\Omega^+_{st}$ and $\Omega^-_{st}$. If this is the case, then by [MY2], there is a
least area annulus in $\BHH$ with boundary $\beta^+_{st}\cup\beta^-_{st}$.

By construction, $\gamma^\pm$ and $\Gamma^\pm_s$ are all pairwise disjoint. By Lemma 2.4, the
geodesic planes $P^\pm$ and the least area planes $\Sigma^\pm_s$ are pairwise disjoint, too. This
implies the boundary of the annulus $\Delta$, $\alpha^+\cup\alpha^-$, is disjoint from $\Sigma^+_s$
and $\Sigma^-_s$. This implies $\Delta\cap\Sigma^\pm_s$ is collection of simple closed curves as
they are least area. Let $\alpha^+_s$ be a simple closed curve in $\Delta\cap\Sigma^+_s$ and
$\alpha^-_s$ be a simple closed curve in $\Delta\cap\Sigma^-_s$. Since $\Sigma^\pm_s$ are least
area planes, the intersection curves $\alpha^+_s$ and $\alpha^-_s$ must be essential curves in
$\Delta$. Otherwise, $\alpha^+_s$ would bound two different disks, one in $\Sigma^+_s$ and the
other one is in $\Delta$. Since $\Sigma^+_s$ and $\Delta$ are both least area already, this cannot
happen by Meeks-Yau exchange roundoff trick \cite{MY2}. Similarly, it is true for $\alpha^-_s$. Let
$\Delta_s\subset\Delta$ be the annulus with boundary $\alpha^+_s\cup\alpha^-_s$.

Now, we make a surgery to get an annulus with boundary $\beta^+_{s0}\cup\beta^-_{s0}$. Let
$D_{\alpha^+_s}$ be the least area disk in $\Sigma^+_s$ with boundary $\alpha^+_s$. Similarly, let
$D_{\alpha^-_s}$ be the least area disk in $\Sigma^-_s$ with boundary $\alpha^-_s$. Since $\Delta$
is least area annulus, the area of $\Delta_s$ is strictly less than the sum of the areas of
$D_{\alpha^+_s}$ and $D_{\alpha^-_s}$. Otherwise, $X= (\Delta-\Delta_s)\cup D_{\alpha^+_s}\cup
D_{\alpha^-_s}$ would be two disks with boundary $\alpha^+\cup\alpha^-$. Moreover, $X$ has singular
circles $\alpha^+_s\cup\alpha^-_s$. By rounding off $X$ along these circles, we get a smaller area
disks, say $X'$. Now, if we put a very thin tube between the disks with a very small area, we get
an annulus $X''$ whose area is less than $\Delta$ with the same boundary. This is a contradiction.
Hence, the area of $\Delta_s$ is strictly less than the sum of the areas of $D_{\alpha^+_s}$ and
$D_{\alpha^-_s}$.

Consider the disks $\Omega^+_{s0}\subset \Sigma^+_{s0}$ with $\partial \Omega^+_{s0} = \beta^+_{s0}$ and $\Omega^-_{s0}\subset
\Sigma^-_{s0}$ with $\partial \Omega^-_{s0} = \beta^-_{s0}$. By construction, $D_{\alpha^\pm_s}$ are subdisks of $\Omega^\pm_{s0}$.  Let
$Y_{s0}= (\Omega^+_{s0}- D_{\alpha^+_s})\cup (\Omega^-_{s0}- D_{\alpha^-_s}) \cup \Delta_s$ be the annulus with boundary
$\beta^+_{s0}\cup\beta^-_{s0}$. Since the area of $\Delta_s$ is strictly less than the sum of the areas of $D_{\alpha^+_s}$ and
$D_{\alpha^-_s}$, the area of $Y_{s0}$ is less than the sum of the areas of the disks $\Omega^+_{s0}$ and $\Omega^-_{s0}$. Then by
\cite{MY2}, there is a least area annulus $A_{s0}$ with $\partial A_{s0}=\beta^+_{s0}\cup\beta^-_{s0}$. Similarly, for any $t\in
[0,\infty)$, there is a least area annulus $A_{st}$ such that $\partial A_{st}=\beta^+_{st}\cup\beta^-_{st}$. Let $\{A_{si}\}$ be a
sequence of least area annuli where $i\in\mathbf{N}$. Notice that the boundary of each annuli in the sequence is
$\beta^+_{si}\cup\beta^-_{si}$ which are simple closed curves in $\Sigma^+_s$ and $\Sigma^-_s$. Since $\PI \Sigma^\pm_s = \Gamma^\pm_s$,
$\beta^+_{si}\rightarrow\Gamma^+_s$ and $\beta^-_{si}\rightarrow\Gamma^-_s$ as $i\rightarrow\infty$. Then, by using the techniques in
\cite{A2}, we can get a subsequence of $\{A_{si}\}$ converging to a complete least area annulus $\mathcal{A}_s$ with $\PI \mathcal{A}_s =
\Gamma^+_s\cup\Gamma^-_s$. Hence, Step 1. follows.
\end{pf}

\begin{figure}[t]
\mbox{\vbox{\epsfbox{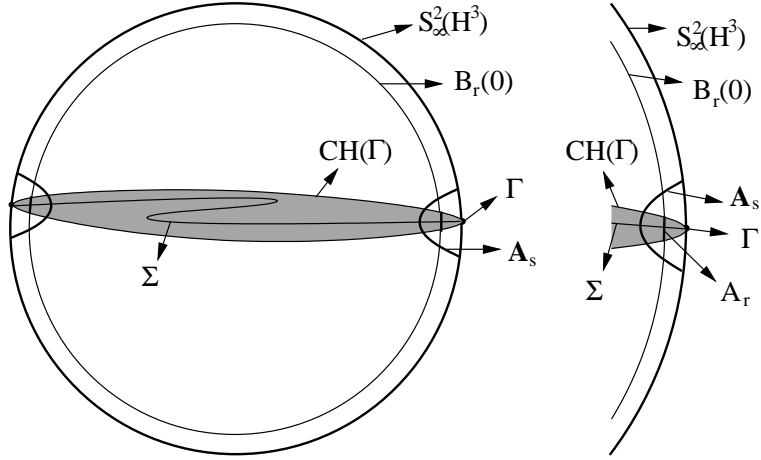}}} \caption{\label{fig:figure1} {The least area annulus $\mathcal{A}_s$ is used as barrier in the proof of
Key Lemma.}}
\end{figure}

\noindent \textbf{Step 2:} ({\em Nonseparating Disks Stays Away from the Center}) There is a function $F$  which is a monotone increasing
function with $F(r)\rightarrow \infty$ as $r \rightarrow \infty$, such that if $D_r$ is a nonseparating disk in $B_r(0)\cap\Sigma$, then
$d(0,D_r)>F(r)$ where $d$ is the distance.\\

\begin{pf}
In the construction in Step 1, we show that for each pair $\Gamma^+_s$ and $\Gamma^-_s$, there is a least area annulus $\mathcal{A}_s$ with
$\PI \mathcal{A}_s = \Gamma^+_s\cup\Gamma^-_s$. Since $\Gamma^\pm_s$ is $C^1$, by \cite{HL} and \cite{To}, $\mathcal{A}_s \cup
\Gamma^+_s\cup\Gamma^-_s$ is a $C^1$ submanifold of the compactification of hyperbolic $3$-space $\overline{\BHH}$. Hence, $\mathcal{A}_s
\cup \Gamma^+_s\cup\Gamma^-_s$ separates $\overline{\BHH}$ into two parts, say $K^+_s$ and $K^-_s$, where $\Gamma\subset K^+_s$. Recall
that $A_r$ is the annulus in $\partial B_r(0)$ such that $int(CH(\Gamma)) \cap \partial B_r(0) \subset A_r$. Define a monotone increasing
function $f: (C',\infty) \rightarrow (C,\infty)$ such that $A_r\subset K^+_{f(r)}$ for $C'$ sufficiently large. Since $\PI \mathcal{A}_s =
\Gamma^+_s\cup\Gamma^-_s \rightarrow \Gamma$ as $s\rightarrow\infty$, the annuli $\{\mathcal{A}_s\}$ escapes to infinity as
$s\rightarrow\infty$. Hence, we can also put the condition $f(r)\rightarrow\infty$ as $r\rightarrow\infty$ on $f$. Now, define a function
$F: (C',\infty) \rightarrow (0,\infty)$ such that $F(r)=d(0,\mathcal{A}_{f(r)})$ where $d$ is the distance in $\BHH$. Clearly, $F$ is a
monotone increasing function and $F(r)\rightarrow\infty$ as $r\rightarrow\infty$.

Now, we claim that if $D_r$ is a nonseparating disk in $B_r(0)\cap\Sigma$, then $d(0,D_r)>F(r)$. To prove this claim, all we need to show
is $D_r\subset K^+_{f(r)}$ for any $r\in (C,\infty)$. In other words, if $A_r$ is in the positive side of $\mathcal{A}_{f(r)}$, i.e.
$A_r\subset K^+_{f(r)}$, then $D_r$ stays in the same side of $\mathcal{A}_{f(r)}$, i.e. $D_r\subset K^+_{f(r)}$ (See Figure 1.). By
assumption $\partial D_r \subset A_r$, and so the boundary of $D_r$ is in the positive side of $\mathcal{A}_{f(r)}$. Assume that $D_r$
intersects $\mathcal{A}_{f(r)}$. Since they are both least area, and the $\partial D_r \cap \mathcal{A}_{f(r)}= \emptyset$, the
intersection is a collection of simple closed curves. Let $\eta$ be such a curve. $\eta$ cannot be essential in $\mathcal{A}_{f(r)}$, since
it bounds a disk in $D_r$, and so $D_r$ will be a separating disk, which contradicts to the assumption. If $\eta$ is not essential in
$\mathcal{A}_{f(r)}$, then this means $\eta$ bounds a disk in $\mathcal{A}_{f(r)}$, too. However, since $\mathcal{A}_{f(r)}$ and $D_r$ are
both least area, this is a contradiction by Meeks-Yau exchange roundoff trick \cite{MY2}. Hence, if $D_r$ is a nonseparating disk in
$B_r(0)\cap\Sigma$, then $d(0,D_r)>F(r)$.
\end{pf}

\end{pf}

\begin{rmk} This lemma is the key point of the main result. Intuitively, this lemma prevents a least area plane to come into the compact part
unnecessarily, where this is very crucial for a plane to be nonproper.
\end{rmk}

\section{Main Result}

In this section, we complete the proof of the main theorem. First, we need a lemma which basically
says that if we have a nonproper least area plane, than we can find arbitrarily large ball such
that the intersection with the least area plane contains infinitely many separating disks.

\begin{lem}
Let $\Sigma$ be an embedded least area plane in $\BHH$ with $\PI \Sigma = \Gamma$ where $\Gamma$ is
a simple closed curve with at least one smooth ($C^1$) point in $\SH$. If $\Sigma$ is not proper,
then for any $R$, there exist $R'>R$ such that the intersection $B_{R'}(0)\cap \Sigma$ contains
infinitely many separating disks.
\end{lem}

\noindent \textbf{Outline:} Assume on the contrary that there exist $R_0>0$ such that there is no
$R>R_0$ such that $B_R(0)\cap\Sigma$ contains infinitely many separating disjoint disks. Fix a
generic $R_1>R_0$. Let $B_{R_1}(0)\cap\Sigma$ contains infinitely many nonseparating disjoint disks
$\{ D_i \}$. Firstly, show that $\{D_i\}$ has an infinite subcollection $\{D_{i_j}\}$ such that
$Area(D_{i_j})>\delta'$ where $\delta'>0$. Then, fix a generic $R_2>R_1$ with $F(R_2)>R_1$. Then
there is a collection of disjoint disks $\{E_l\}$ in $B_{R_2}(0) \cap \Sigma$ such that for any
$i_j$ there is an $l$ with $D_{i_j} \subset E_l$. By using the area bound, show that $\{E_l\}$ is
an infinite collection of disjoint disks, and by the assumption, for all but finitely many, they
are nonseparating. Let $E_{l_1}$ be a such nonseparating disk. Since $D_{i_k}\subset E_{l_1}$,
$d(0,E_{l_1})< R_1$. Since $E_{l_1}$ is nonseparating, and by Key Lemma, $d(0,E_{l_1})>
F(R_2)>R_1$. This is a contradiction.\\

\begin{pf}
By Lemma 3.1, there exist $r_0>0$ such that for a generic $r>r_0$, $B_r(0)\cap\Sigma$ contains infinitely many disjoint disks. Assume that
there exist $R_0> r_0$ such that there is no $R>R_0$ such that $B_R(0)\cap\Sigma$ contains infinitely many separating disjoint disks.

Let $R_1>R_0$ be a generic radius in the sense of Lemma 3.1. i.e. $B_{R_1}(0)\cap\Sigma$ contains infinitely many nonseparating disjoint
disks $\{ D_i \}$. Let $\gamma_i = \partial D_i$ be the pairwise disjoint simple closed curves in the annulus $A_{R_1} \subset \partial
B_{R_1}(0)$. Since $\{D_i\}$ are nonseparating, $\{\gamma_i\}$ are not essential in $A_{R_1}$. Let $\Omega_i$ be the disk in $A_{R_1}$ with
$\partial \Omega_i = \gamma_i$. We claim that $\{\gamma_i\}$ has an infinite subsequence $\{\gamma_{i_j}\}$ with
$\Omega_{i_j}\supset\Omega_{i_k}$ for any $i_j<i_k$.

Assume on the contrary that there is no such subsequence. Since $\{\gamma_i\}$ is a collection of pairwise disjoint curves in $A_{R_1}$,
$\{\gamma_i\}$ must have an infinite subsequence $\{\gamma_{i_j}\}$ with $\Omega_{i_j}\cap\Omega_{i_k}=\emptyset$ for any $i_j\neq i_k$. We
can also assume that any curve in the sequence is an outermost curve, i.e. $\Omega_{i_j}\nsubseteq\Omega_k$ for $k\neq i_j$. Since the area
of $A_{R_1}$ is finite, $\sum_{j=1}^\infty Area(\Omega_{i_j})< \infty$. This implies as $i_j \rightarrow \infty$,
$Area(\Omega_{i_j})\rightarrow 0$. Since $\{D_{i_j}\}$ are least area, $Area(D_{i_j})\rightarrow 0$ as well. Let $\epsilon>0$ be a
sufficiently small number with $R_1-\epsilon$ is also a generic in the sense of Lemma 3.1. As $R_1-\epsilon>R_0$, by assumption,
$B_{R_1-\epsilon}(0)\cap\Sigma$ must contain infinitely many nonseparating disjoint disks $\{ E_i \}$. However, as $i_j \rightarrow
\infty$, $Area(\Omega_{i_j})\rightarrow 0$ and $\{\gamma_{i_j}\}$ are outermost curve, $B_{R_1-\epsilon}(0)\cap\Sigma$ can have only
finitely many disks. This is a contradiction. Hence, we can assume that $\{\gamma_i\}$ has a subsequence $\{\gamma_{i_j}\}$ with
$\Omega_{i_j}\supset\Omega_{i_k}$ for any $i_j<i_k$.

By using the similar ideas in previous paragraph, we can also assume that there is a disk $\Omega_\infty$ in $A_{R_1}$ with
$Area(\Omega_\infty) > \delta$ and for any $i_j$, $\Omega_{i_j}\supset\Omega\infty$. This is because if this is not the case, then we will
be in the situation as $i_j \rightarrow \infty$, $Area(\Omega_{i_j})\rightarrow 0$ as before, and similar ideas will give us a
contradiction. Hence, we can assume that for any $i_j$, $Area(\Omega_{i_j})>\delta$. Clearly, there is a $\delta'>0$ such that
$Area(D_{i_j})>\delta'$.

Now, let $R_2>R_1$ be also generic in the sense of Lemma 3.1 and $F(R_2)>R_1$. Since $B_{R_1}(0)
\subset B_{R_2}(0)$, the $B_{R_1}(0) \cap \Sigma \subset B_{R_2}(0) \cap \Sigma$. Hence, there is a
collection of disjoint disks $\{E_l\}$ in $B_{R_2}(0) \cap \Sigma$ such that for any $i_j$ there is
an $l$ with $D_{i_j} \subset E_l$. We claim that the collection $\{E_l\}$ contains infinitely many
disjoint disks.

If the collection $\{E_l\}$ has only finitely many disjoint disks, then there is an $l_0$ such that $E_{l_0}$ contains infinitely many
disks in $\{D_{i_j}\}$. By the proof of Lemma 3.1, for any $l$, the disk $E_l$ must have finite area. Since, for any $i_j$,
$Area(D_{i_j})>\delta'$ and the area of $E_{l_0}$ is finite, this is a contradiction. Therefore, the collection $\{E_l\}$ contains
infinitely many disjoint disks.

Since $R_2>R_0$, all but finitely many disks in the collection $\{E_l\}$ must be nonseparating. Let
$E_{l_1}$ be a nonseparating disk in the collection. Let $D_{i_k}\subset E_{l_1}$. Since $D_{i_k}
\subset B_{R_1}(0)\cap\Sigma$, $d(0,D_{i_k})< R_1$. Hence, $d(0,E_{l_1})< R_1$. However, $E_{l_1}$
is nonseparating, and by Key Lemma, $d(0,E_{l_1})> F(R_2)$. Since $F(R_2)>R_1$, this is a
contradiction. Hence, the proof follows.

\end{pf}

Now, we can prove the main result of the paper.

\begin{thm}
Let $\Sigma$ be a complete embedded least area plane in $\BHH$ with $\PI\Sigma=\Gamma$ where
$\Gamma$ is a simple closed curve in $\SH$ with at least one smooth ($C^1$) point. Then, $\Sigma$
must be proper.
\end{thm}

\begin{pf}
Assume that $\Sigma$ is not proper. Then by Lemma 5.1, for any $R>0$, there exist $R'>R$ such that the intersection $B_{R'}(0)\cap \Sigma$
contains infinitely many separating disks. Let $\beta$ be a path from $\partial^+ CH(\Gamma)$ to $\partial^- CH(\Gamma)$ through $0$. Let
$l$ be the length of $\beta$. Let $R_1>0$ be so that $F(R_1)>l$, and $B_{R_1}(0)\cap \Sigma$ contains infinitely many pairwise disjoint
separating disks $\{D_i\}$. Let $x_i\in D_i$ for any $i$. Let $\{\gamma_{ij}\}$ be the family of paths in $\Sigma$ between $x_i$ and $x_j$.
Let $R_2= \inf_{R>R_1} \{ R \ | \ \exists i,j>0, \exists \gamma_{ij}\subset\Sigma, \gamma_{ij} \subset B_R(0) \}$. In other words,
$B_{R_2}(0)$ is the smallest closed ball among the balls $B_R(0)$ where at least two of the disjoint disks in $B_{R_1}(0)\cap \Sigma$ can
be connected in $B_R(0)\cap \Sigma$. Say, we can connect $D_i$ and $D_j$ in $B_{R_2}(0)$ via path $\gamma_{ij}$. Let $E$ be the component
of $B_{R_2}(0)\cap \Sigma$ containing $D_i$ and $D_j$, i.e. $D_i\cup D_j \subset E \subset B_{R_2}(0)\cap \Sigma$. Because of the
assumption on $R_2$, $\gamma_{ij}\cap\partial B_{R_2}(0)\neq \emptyset$. Hence, $\partial E$ in $\partial B_{R_2}(0)$ is a nonsimple closed
curve (with degenerate point $\gamma_{ij}\cap\partial B_{R_2}(0)$). Then, by changing the center a little bit from $0$ to $0'$, if
necessary; we can find a sufficiently small $\epsilon$ such that $R_2'=R_2+\epsilon$ is a generic radius in the sense of Lemma 3.1, and the
component $E'$ in $B_{R_2'}(0)$ containing $D_i$ and $D_j$ is a nonseparating disk. Since $D_i$ is a separating disk in $B_{R_1}(0)$,
$d(0,D_i)<l$. This implies $d(0,E')<l$. However, since $E'$ is a nonseparating disk in $B_{R_2'}(0)$, $d(0,E')>F(R_2')$. Since
$F(R_2')>F(R_1)>l$, this is a contradiction. The proof follows.

\end{pf}

\section{Concluding Remarks}

As it is mentioned in the introduction, even though many experts of the field has questioned the problem, there is a very few results about
the properly embeddedness of least area planes in $\BHH$. For example, it is still not known if there exists a nonproperly embedded least
area plane in $\BHH$ whose asymptotic boundary is a simple closed curve in $\SH$. On the other hand, there is a construction by Freedman
and He communicated to Gabai of a nonproper least area plane in $\BH^3$.  It is not clear to the author as to how a plane constructed in
this manner can have limit set a simple closed curve.

Recently, Colding and Minicozzi solved an analogous question in $\BR^3$ in \cite{CM}. They proved
Calabi-Yau Conjectures for embedded surfaces by relating intrinsic distances and extrinsic
distances of the minimal disk. As a corollary, they proved that any complete embedded minimal plane
in $\BR^3$ must be proper. Our approach is very different from them as they use purely analytic
methods, while our techniques are purely topological. When starting this problem, our aim is to
prove the following conjecture.\\

\noindent \textbf{Conjecture:} Let $\Sigma$ be a complete embedded least area plane in $\BHH$ with
$\PI\Sigma=\Gamma$ where $\Gamma$ is a simple closed curve in $\SH$. Then, $\Sigma$ must be
proper.\\

In this paper, we proved this statement with the existence of a smooth point condition on $\Gamma$.
We needed this condition for the Key Lemma to show existence of a least area annulus linking $\PI
\Sigma$. If one can bypass this without the smooth point condition, then the theorem can be proved
in full generality. On the other hand, since the smooth point condition also means the finite
thickness of $CH(\Gamma)$ in one direction, this might be an essential point for the result.


\begin{thebibliography}{MSY}


\bibitem[A1]{A1} M. Anderson, {\em Complete minimal varieties in hyperbolic space}, Invent. Math. {\bf 69}, 477--494 (1982).

\bibitem[A2]{A2} M. Anderson, {\em Complete minimal hypersurfaces in hyperbolic n-manifolds}, Comment. Math. Helv. {\bf 58}, 264--290
(1983).

\bibitem[CM]{CM} T.H. Colding and W.P. Minicozzi, {\em The Calabi-Yau conjectures for embedded
surfaces}, to appear in Annals of Math.

\bibitem[Co1]{Co1} B. Coskunuzer, {\em Uniform 1-cochains and Genuine Laminations}, Topology {\bf 45} (2006) 751--784.

\bibitem[Co2]{Co2} B. Coskunuzer, {\em Generic Uniqueness of Least Area Planes in Hyperbolic Space}, Geom. \& Topology
{\bf 10} (2006) 401--412.

\bibitem[Co3]{Co3} B. Coskunuzer, {\em Properly Embedded Least Area Planes in Gromov Hyperbolic 3-Spaces},
to appear in Proc. Amer. Math. Soc.

\bibitem[EM]{EM}D.B.A. Epstein and A. Marden, {\em Convex hulls in hyperbolic space, a theorem of Sullivan, and measured pleated
surfaces},  London Math. Soc. Lecture Note Ser. {\bf 111}, (1987) 113--253.

\bibitem[Ha]{Ha} J. Hass, {\em Intersections of Least Area Surfaces}, Pacific J. Math. {\bf 152} (1992) 119--123.

\bibitem[HL]{HL} R. Hardt and F.H. Lin, {\em Regularity at infinity for absolutely area minimizing hypersurfaces in hyperbolic space}, Invent.
Math. {\bf 88} (1987) 217--224.

\bibitem[Ga]{Ga} D. Gabai, {\em On the geometric and topological rigidity of hyperbolic $3$-manifolds}, J. Amer. Math. Soc. {\bf 10} (1997) 37--74.

\bibitem[MY1]{MY1} W. Meeks and S.T. Yau, {\em Topology of three-dimensional manifolds and the embedding
problems in minimal surface theory}, Ann. of Math.  {\bf 112} (1980) 441--484.

\bibitem[MY2]{MY2} W. Meeks and S.T. Yau, {\em The classical Plateau problem and the topology of three manifolds}, Topology {\bf 21}
(1982) 409--442.

\bibitem[So1]{So1} T. Soma, {\em Existence of least area planes in hyperbolic 3-space with co-compact metric}, Topology {\bf 43} (2004) 705--716.

\bibitem[So2]{So2} T. Soma, {\em Least area planes in Gromov hyperbolic 3-spaces with co-compact metric}, Geom. Dedicata {\bf 112} (2005) 123--128.

\bibitem [To]{To} Y. Tonegawa, {\em Existence and regularity of constant mean curvature
hypersurfaces in hyperbolic space}, Math. Z. 221, (1996) 591--615.

\end{thebibliography}
\end{document}